\documentclass[12pt]{amsart}
\usepackage{setspace}
\usepackage{amsmath,amsthm,amssymb,graphicx,mathtools,tikz,hyperref}
\usepackage{ dsfont }
\usetikzlibrary{positioning}
\newcommand{\n}{\mathbb{N}}
\newcommand{\z}{\mathbb{Z}}

\newcommand{\real}{\mathbb{R}}

\newcommand{\p}{\mathbb{P}}
\newcommand{\e}{\mathbb{E}}

\newtheorem{theorem}{Theorem}[section]
\newtheorem{proposition}[theorem]{Proposition}

\newtheorem{lemma}[theorem]{Lemma}
\newtheorem{conjecture}[theorem]{Conjecture}

\usepackage{mdframed}
\usepackage{lipsum}
\newmdtheoremenv{defi}{Definition}
\newmdtheoremenv{theo}{Theorem}
\newmdtheoremenv{prop}{Proposition}
\newmdtheoremenv{lemm}{Lemma}
\newmdtheoremenv{coro}{Corollary}
\newmdtheoremenv{conj}{Conjecture}
\newmdtheoremenv{constr}{Construction}

\theoremstyle{definition}
\newtheorem*{definition}{Definition}

\usepackage{times}

\title{Lower bounds for the universal TSP on the plane}
\author{Cosmas Kravaris}
\email{ck6221@princeton.edu}
\address{Mathematics Department, Princeton University, Fine Hall, Washington Road, Princeton, NJ 08544-1000, USA}

\begin{document}
\maketitle

\begin{abstract}
We show a lower bound for the universal traveling salesman heuristic on the plane:
for any linear order on the unit square $[0,1]^2$,
there are finite subsets $S \subset [0,1]^2$ of arbitrarily large size
such that the path visiting each element of $S$ according to the linear order has length $\geq C \sqrt{\log |S| / \log \log |S|}$ times the length of the shortest path visiting each element in $S$. ($C>0$ is a constant that depends only on the linear order.)
This improves the previous lower bound $\geq C \sqrt[6]{\log |S| / \log \log |S|}$ of Hajiaghayi, Kleinberg and Leighton (SODA 2006).
The proof establishes a dichotomy about any long walk on a cycle: the walk either zig-zags between two far away points, or else for a large amount of time it stays inside a set of small diameter.
\end{abstract}

\section{Introduction.}

Suppose you have a list of items to buy from the grocery store and want to pick them up as quickly as possible.
So you want to order all the items in the list to minimize the sum of all the distances from the first item to the second, from the second item to the third, and so forth.
This is the famous \textbf{traveling salesman problem} in combinatorial optimization. 
Its computational complexity has been extensively studied; see, for instance, the book \cite{korte2008combinatorial}, the paper \cite{karlin2021slightly}, and references within.
In this paper, we study the competitive ratio of a specific, natural, and simple heuristic for the traveling salesman problem.

If you visit the same grocery store every week, then you might want to do some \textbf{pre-processing} before you know the list of items to buy. One strategy is to place a universal (or global) linear order on the entire space, such as a tour that visits every single aisle. Given any list of items, you may sort them according to this universal linear order (i.e. write the items on your phone according to the order you would visit them if you took the entire tour). Then visit the items according to this order (going directly from each item to the next, skipping any aisles where there are no items). This heuristic is a very natural algorithm for the traveling salesman problem, and many people follow it instinctively whenever they shop regularly at the same grocery store. Surprisingly, we still do not fully understand the competitive ratio analysis for this simple heuristic, even on the plane.

Here is the precise formulation of the \textbf{universal traveling salesman problem}.
It was introduced by Bartholdi and Platzman \cite{platzman1989spacefilling} who studied it on the plane, and was first studied on general metric spaces in \cite{jia2005universal}.
Let $(\mathcal{M},d)$ be a metric space.
For any linear order $\leq$ on $\mathcal{M}$ and any finite subset $\{s_1,\dots s_n\}\subset \mathcal{M}$, indexed such that $s_1 < s_2 < \dots  < s_n$, the \textbf{cost} of visiting each point in the set according to $\leq$ is
$$cost_{\leq}(\{s_1,\dots s_n\}) := \sum_{i=1}^{n-1} d(s_{i},s_{i+1}).$$
The \textbf{traveling salesman problem (TSP)} asks for the smallest cost among all possible ways we can order the points $S = \{s_1,\dots ,s_n\}$.
We write
$$tsp(\{s_1,\dots ,s_n\}) := \inf_{\pi \in Sym_n} \sum_{i=1}^{n-1} d(s_{\pi(i)},s_{\pi(i+1)}),$$
where $Sym_n$ denotes the set of all permutations on $\{1,\dots ,n\}$.

In this paper, we follow Bartholdi-Platzman \cite{platzman1989spacefilling} and \cite{hajiaghayi2006improved, bhalgat2011optimal, erschler2023spaces, erschler2023assouad, mitrofanov2022total} by defining the \textbf{TSP competitive ratio function} (or \textbf{order ratio function}) of a linear order $\leq$ on $\mathcal{M}$ as
$$\rho_n(\mathcal{M},\leq) := \sup_{S \subset \mathcal{M}:\;|S|\leq n} \dfrac{cost_{\leq}(S)}{tsp(S)},\;\;\;for\;all\;n \in \n,$$
i.e. we study the worst competitive ratio \textit{as a function of the size of the input}.
Later literature also studies another way to parametrize the competitive ratio by the size of the entire metric space $\mathcal{M}$ - see subsection \ref{intro:variants} for a discussion on this variant.

In \cite{platzman1989spacefilling}, Bartholdi and Platzman used
the Sierpinski space-filling curve $f: [0,1] \to [0,1]^2$ to define a linear order $\leq_{BP}$ on $[0,1]^2$ by
$$s \leq_{BP} s' \iff \min\{t\in [0,1]: f(t) = s\} \leq \min\{t\in [0,1]: f(t) = s'\}.$$
They showed that (with respect to the standard Euclidean metric on $[0,1]^2$),
$$\rho_{n}([0,1]^2,\leq_{BP}) \lesssim \log n\;\;\;for\;all\;n \in \n,$$
and conjectured that $\rho_n([0,1]^2,\leq_{BP}) \lesssim 1$.
\footnote{Recall that for two sequences $\{a_n\}_n, \{b_n\}_n \subset \real^+$ one has $a_n \lesssim b_n$ if and only if $b_n \gtrsim a_n$ if and only if there exists $0<C<\infty$ such that $a_n \leq C b_n$ for all $n \in \n$. Also, we write $a_n \asymp b_n$ when $a_n \lesssim b_n$ and $b_n \lesssim a_n$.}
Their conjecture was disproved by Bertsimas and Grigni \cite{bertsimas1989worst} who showed that 
$$\rho_n([0,1]^2,\leq_{BP}) \gtrsim \log n\;\;\;for\;all\;n \in \n.$$
They conjectured that their bound is tight among all possible linear orders.
In terms of the function $\rho_n([0,1]^2,\leq)$ studied in \cite{platzman1989spacefilling} and \cite{hajiaghayi2006improved} the conjecture is as follows (see subsection \ref{intro:variants} for a variant):

\begin{conjecture}[Bertsimas, Grigni]\label{conjecture:relative}
    For any linear order $\leq$ on the unit square $[0,1]^2$
    $$\rho_n([0,1]^2,\leq) \gtrsim \log n.$$
\end{conjecture}

In \cite{hajiaghayi2006improved}, Hajiaghayi, Kleinberg and Leighton showed that any linear order $\leq$ on $[0,1]^2$ has
$$\rho_n([0,1]^2,\leq) \gtrsim \sqrt[6]{\dfrac{\log n}{\log \log n}}.$$
In \cite{eades2020optimal}, Eades and Mestre proved Conjecture \ref{conjecture:relative} for a special family of orders called "hierarchical". 

The purpose of this paper is to prove the following theorem:
\begin{theorem}\label{main:theorem}
    For any linear order $\leq$ on the unit square $[0,1]^2$,
    $$\rho_n([0,1]^2,\leq) \gtrsim \sqrt{\dfrac{\log n}{\log \log n}},$$
    that is, for any $n \in \n$ there exists an $n$-point set $S \subset [0,1]^2$ with
    $$cost_{\leq}(S) \geq C \sqrt{\dfrac{\log n}{\log \log n}}\; tsp(S),$$
    where $C>0$ is a constant independent of $n$ and $S$.
\end{theorem}

Some remarks: First, the lower bound does not improve if we replace the unit square $[0,1]^2$ with the unit cube $[0,1]^d$ in higher dimension $d \geq 3$. 
It would also be interesting to understand the competitive ratio of metric spaces of fractal dimension between $1$ and $2$.

Secondly, there is no necessity to work with the infinite metric space $[0,1]^2$, even though this is how the literature frames the question.
It suffices to show that for each $n \in \n$ there exists $D \in \n$ such that any linear order $\leq$ on the $D \times D$ grid contains a subset of size $n$ with competitive ratio $\gtrsim \sqrt{\log n/ \log \log n}$.
(The converse is also true as observed in \cite{erschler2023assouad}. By the compactness theorem in logic, if there exists an efficient order on every finite subset of a metric space, then there exists an efficient order on the entire space.)

Finally, see subsection \ref{Remark:Barrier_of_1_over_2} for a heuristic explanation that $1/2$ is the barrier exponent in our proof of Theorem \ref{main:theorem}. 

\subsection{Related literature on the universal TSP}\label{subSect:related_literature}

The universal traveling salesman problem has been studied on a plethora of spaces in metric geometry.
There are upper bounds (i.e. constructions of efficient orders) for the plane \cite{platzman1989spacefilling, christodoulou2017improved}, doubling metric spaces \cite{jia2005universal}, spaces of finite Assouad-Nagata dimension \cite{erschler2023assouad}, hyperbolic graphs \cite{erschler2023spaces}, minor-excluded graphs \cite{hajiaghayi2006improved}, and general $n$-point metrics \cite{jia2005universal}.
There are lower bounds (i.e. impossibility of efficient orders) for the plane \cite{bertsimas1989worst, hajiaghayi2006improved, eades2020optimal}, Ramanujan graphs of high girth \cite{bhalgat2011optimal, gorodezky2010improved}, expanders \cite{erschler2023assouad} and higher dimensional spheres and their variants \cite{erschler2023assouad}.

In \cite{erschler2023assouad}, Erschler and Mitrophanov showed a structural dichotomy about all linear orders on doubling metric spaces (such as the plane): every linear order $\leq$ either has $\rho_n([0,1]^2,\leq) \lesssim \log n$ or else $\rho_n([0,1]^2,\leq)=n$ for all $n\in \n$ (so there is an \textbf{order gap}).
We use this result in the proof of Theorem \ref{main:theorem}, otherwise, we only get the lower bound for infinitely many values of $n \in \n$ instead of all values of $n$.

Other structural properties that are studied in the literature are the order breakpoint introduced in Erschler-Mitrophanov \cite{erschler2023spaces} (this is the smallest $n \in \n$ such that $\rho_n(\mathcal{M},\leq) < n$ for all linear orders $\leq$) 
and Mitrophanov's characterization of finite covering dimension via snakes \cite{mitrofanov2022total}.

The \textit{universal heuristic} has also been studied for other combinatorial optimization problems, such as the Steiner tree problem and the set cover problem; see \cite{bhalgat2011optimal,jia2005universal}.

\subsection{Variant on the universal TSP competitive ratio}\label{intro:variants}

The works \cite{bertsimas1989worst, christodoulou2017improved, eades2020optimal, bhalgat2011optimal, gorodezky2010improved, jia2005universal} study a global variant of the competitive ratio, namely the quantity $\sup_{S \subset \mathcal{M}} cost_{\leq}(S)/tsp(S)$ as a function on the number of points $|\mathcal{M}|$ of the metric space $\mathcal{M}$.

The variant of Conjecture \ref{conjecture:relative} in terms of the above parameterization states that for any $n \in \n$ and any linear order $\leq$ on the $n \times n$ grid $[n]^2$ there exists a subset $S \subset [n]^2$ with competitive ratio $cost_{\leq}(S)/tsp(S) \gtrsim \log n$.
This version of the conjecture is clearly stronger than Conjecture \ref{conjecture:relative}.
It was disproved by Christodoulou and Sgouritsa \cite{christodoulou2017improved} who showed that for each $n \in \n$ there exists a linear order $\leq$ on $[n]^2$ with global competitive ratio $\sup_{S \subset \mathcal{M}} cost_{\leq}(S)/tsp(S) \lesssim \log n/\log \log n$.

The argument of \cite{hajiaghayi2006improved} gives that for each linear order $\leq$ on  $[n]^2$ there exists a subset  $S \subset [n]^2$ with competitive ratio $cost_{\leq}(S)/tsp(S) \gtrsim \sqrt[6]{\log n/\log \log n}$.
The proof of Theorem \ref{main:theorem} does not imply any improved result on the above variation; the gap between $\Omega (\log n/ \log \log n)^{1/6}$ and $O(\log n/ \log \log n)$ persists for the global competitive ratio of linear orders on the $n \times n$ grid.

\subsection{Discussion on the proof and the new insight}\label{intro:proof_discussion}

Similar to \cite{hajiaghayi2006improved},
there are two types of obstructions to a small competitive ratio function, backtracks and zig-zags, and we show that at least one of the two must occur.
Each obstruction gives a set $S$ with $cost_{\leq}(S) \geq C \sqrt{\log |S|/ \log \log |S|}\; tsp(S)$. (For zig-zags, the competitive ratio is worse as a function of $|S|$.)

The first type is a \textbf{backtracking set}: all points are close to a line $L$ and the ordering of $S$ backtracks a lot for many relevant distance scales. (The idea goes back to \cite{bertsimas1989worst}.)
Each dyadic square in $[0,1]^2$ contains an individual backtrack in some direction,
and the set $S$ is the union of backtracks which are close to $L$ and have the same direction as $L$.
The new insight on this side of the dichotomy is the definition of a "backtrack" which allows us to improve the lower bound.
See subsection \ref{Remark:comparingDefinitions} for a comparison with the definition used in \cite{hajiaghayi2006improved}. 
Given a number of equi-distributed angles $M\in \n$, a length $0<l<1$ and a width $0<w<l$, 
a \textbf{backtrack} for $\leq$ in $[0,1]^2$ consists of two rectangles $R_1$ and $R_2$ of length $l$ and width $w$ which lie inside the $w$-neighborhood of a line $L$ of angle $\in\{2\pi/M,4\pi/M,6\pi/M,\dots 2\pi\}$, 
and a point $p \in [0,1]^2$ which has distance $w/2 < d(p,L) < w$ to the line $L$ and lies between $R_1$ and $R_2$ such that
$$p<q \;\;\;for\;all\;q \in R_1 \cup R_2.$$
This definition is stronger than the one given in \cite{hajiaghayi2006improved} (which is insufficient for our purposes; read two paragraphs ahead),
and the analysis of the competitive ratio for the backtracking set follows \cite{hajiaghayi2006improved}.

The second type of obstruction is a \textbf{zig-zag}: all points in $S$ lie close to one of two line segments, and the linear order on $S$ jumps back and forth many times between the two segments.
The construction of the zig-zag set is new and short.
Suppose that $[0,1]^2$ (or some dyadic subsquare) has no backtrack.
Starting from the point $p:=[3/4,1/2]$ we iteratively use the fact that $p$, the radial line $L$ at $p$ with respect to the center $(1/2,1/2)$, and the two rectangles $R_1$ and $R_2$ on $L$ of distance $1/M$ from $p$ do not form a backtrack.
This means that there exists $q \in R_1 \cup R_2$ such that $p>q$ and now apply the same argument to $q$.
We obtain a \textbf{spiral chain} of points $p_1 > p_2 > p_3 > \dots .$, each of distance $\asymp 1/M$ from the previous one.
The points move further and further from the center $(1/2,1/2)$ of $[0,1]^2$.
The uniform smoothness of the space $l_2^2$ means that this chain will have $\gtrsim M^2$ points before exiting the square.

The \textbf{key insight
is a new dichotomy} about walks of length $M^2$ on the $M$-cycle.
For any $1 \leq s \leq M^{1/3}$, either there are two points of distance $s^2$ apart such that the walk zig-zags between them $M/2s$ many times, 
or else for $s^3$ consecutive steps, the walk is confined within a set of diameter $6 s^2 + 2$. (The parameter dependencies are tight up to constant factors.)
The following two examples illustrate the dichotomy.
At one extreme, the \textbf{constant walk} always moves clockwise on the $M$-cycle and makes $M$ total rotations.
This walk zig-zags $M$-times between the distant points $1$ and $\lfloor M/2\rfloor$.
At the other extreme, the \textbf{shy walk} alternates between moving clockwise and counterclockwise, oscillating at two neighboring points for the entire walk.
This entire walk is a confined inside a set of diameter $1$.
The spiral chain  $p_1 > p_2 > p_3 > \dots .$ corresponds to a walk of length $M^2$ on the $M$-cycle, and each of the two cases in the dichotomy gives us a zig-zag between two radial line segments (either two long line segments which are far away, or else two small line segments which are close to each other).

\section{Overview of the proof}\label{Sec:ProofOverview}

Throughout the paper $([0,1]^2,d)$ is the unit square with the standard Euclidean metric $d((a_1,a_2), (b_1,b_2)) := \sqrt{(a_1-b_1)^2 + (a_2-b_2)^2}$ for every $0 \leq a_1,a_2,b_1,b_2 \leq 1$.

Fix a linear order $\leq$ on $[0,1]^2$. First note that
\\if we show $\rho_{2^r}([0,1]^2,\leq) \gtrsim \sqrt{\dfrac{r}{\log r}}$ for all $r \in \n$, 
\\then we get $\rho_n([0,1]^2,\leq) \gtrsim \sqrt{\dfrac{\log n}{\log \log n}}$ for all $n \in \n$.
\\We argue by contradiction.
Assume this is not the case.
Then there exists a subsequence $(r^{(k)})_{k \in \n}$ of $(r)_{r \in \n}$ such that
$$\rho_{2^{r^{(k)}}}([0,1]^2,\leq)= o\left(\sqrt{\dfrac{r^{(k)}}{\log r^{(k)}}}\right).$$
This implies that there exists a sequence $(M_k)_k \subset \n$ with $\lim_{k \to \infty} M_k = \infty$ such that
$$\rho_{2^{r^{(k)}+2}}([0,1]^2,\leq) < 10^{-4} \sqrt{\dfrac{r^{(k)}}{M_k^9 \log r^{(k)}}}.$$
We will prove the following dichotomy:
\begin{proposition}[The dichotomy which will imply the main theorem]\label{main-prop}.\\
For any linear order $\leq$ on $[0,1]^2$, any $r \in \n$, and any $M \in \n$ 
such that $180^2 < M \leq 10^{-5} (r/\log r)^{1/9}$,
\\\textbf{CASE A}: either there exists a subset $S \subset [0,1]^2$ of size $|S| \leq 2^{r+2}$ with 
$$\dfrac{cost_{\leq}(S)}{tsp(S)}\geq 10^{-4} \sqrt{\dfrac{r}{M^9 \log r}},$$
\\\textbf{CASE B}: or else there exists a subset $S \subset [0,1]^2$ of size $|S| \lesssim M^{2/3}$ with 
$$\dfrac{cost_{\leq}(S)}{tsp(S)} \gtrsim M^{1/3} \gtrsim \sqrt{|S|}.$$
\end{proposition}

Let's see how the above Proposition \ref{main-prop} gives us Theorem \ref{main:theorem}.
We can modify our sequence $M_k$ to be $\min\{M_k, 10^{-5} (r^{(k)}/\log r^{(k)})^{1/9}\}$ and take sufficiently large $k$ such that $M_k > 180^2$ so that the hypothesis of Proposition \ref{main-prop} is satisfied.
Since the first part of the dichotomy cannot hold, there are infinitely many values $n \in \n$, namely $n=M_{k}^{2/3}$, such that
$$\rho_{n}([0,1]^2,\leq) \gtrsim \sqrt{n} \gg \sqrt{\dfrac{\log n}{\log \log n}}.$$
In particular, we obtain Theorem \ref{main:theorem} for infinitely many values of $n$.

If we use the order gap property for doubling metric spaces due to Erschler and Mitrophanov \cite{erschler2023assouad}, we obtain the lower bound for all values of $n \in \n$:
\begin{theorem}[Erschler, Mitrophanov]
    .\\For any linear order $\leq$ on a doubling metric space $(\mathcal{M},d)$ (such as $([0,1]^2,d)$) we have:
    \\\textbf{either} $\rho_{n}(\mathcal{M},\leq) \lesssim \log n$ for all $n \in \n$
    \\\textbf{or else} $\rho_{n}([\mathcal{M},\leq) = n$ for all $n \in \n$.
\end{theorem}
The bound
$\limsup_{n \to \infty} \rho_{n}([0,1]^2,\leq)/\sqrt{n} > 0$
means that the first part of the dichotomy cannot hold.
This forces $\rho_{n}([0,1]^2,\leq) = n \gg \sqrt{\dfrac{\log n}{\log \log n}}$ for all $n \in \n$, arriving at our contradiction.

In the rest of the paper, we show Proposition \ref{main-prop}.
Fix the linear order $\leq$, and the two parameters: 
\\$r \in \n$, the \textbf{number of scales}, and $M \in \n$, the \textbf{number of angles} which is sufficiently large: $M> 180^2$. 
(We will need this assumption for the spiral chain construction in Lemma \ref{spiral:chain:construction}.)
\\At the very end of the proof (Section 7), we use the assumption 
$M \leq 10^{-5} (r/\log r)^{1/9}$.

\begin{definition}\label{main:definition}
A \textbf{backtrack} $(p,L,\sigma,R_1,R_2)$ of length $l>0$ and width $w>0$ consists of:
\\a line $L$ of angle $\{\dfrac{2\pi}{M}, 2\dfrac{2\pi}{M}, 3\dfrac{2\pi}{M},\dots , 2\pi\}$,
\\a point $p \in [0,1]^2$ of distance $d(p,L) < w/4$ from $L$
\\a strip $\sigma$ which is the $w/2$-neighborhood of line $L$,
\\and two rectangular regions $R_1, R_2 \subset [0,1]^2$ of width $w$ and length $l$ which lie in $\sigma$ in opposite sides of $p$ satisfying:
$$ p < q \;\;\; for \; all \; q \in R_1 \cup R_2.$$  
\end{definition}

Fix two parameters, the \textbf{length} $0<l<1/2$ and the \textbf{width} $0<w<l$ that depend on $r$ and $M$ (at the very end of the proof we will optimize and take $l \asymp M^{-4}$ and $w \asymp \sqrt{M \log r /r}$).

For each scale $t=0,\dots ,r$ consider all dyadic subsquares of $[0,1]^2$ of scale $t$.
There are $2^{2t}$ many squares $Q$ in total and they partition $[0,1]^2$.
Each square is of the form $Q = [a_1,a_1+2^{-t}] \times [a_2,a_2+2^{-t}]$ for some $a_1,a_2 \in \{0,2^{-t},2 \times 2^{-t},3 \times 2^{-t},\dots ,(2^t-1)\times 2^{-t}\}$.

We say a dyadic square $Q$ of scale $t$ \textbf{contains a backtrack} $(p,L,\sigma,R_1,R_2)$ when $\{p\}\cup R_1 \cup R_2 \subset Q$ and
$(p,L,\sigma,R_1,R_2)$ is a backtrack of width $2^{-t}w$ and length $2^{-t}l$.
(Note that the width and length are scaled according to the sidelength of $Q$.)
The proof now splits into the following two cases, corresponding to the same cases A and B in Proposition \ref{main:definition}.
\\ \textbf{CASE A: every dyadic square contains a backtrack}. 
\\In this case, we get a backtracking set. We analyze it later.
\\ \textbf{CASE B: there exists a dyadic square with no backtrack}.
\\In this case, we get a zig-zag. Let's begin with this case.

\begin{figure}[h]
	\begin{center}
		\includegraphics[scale=0.20]{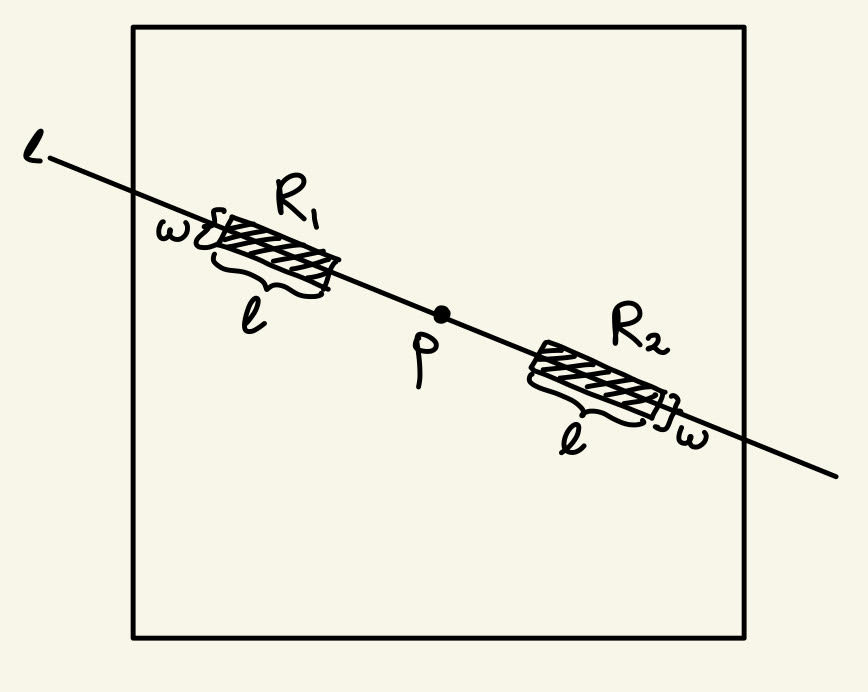}
        \includegraphics[scale=0.17]{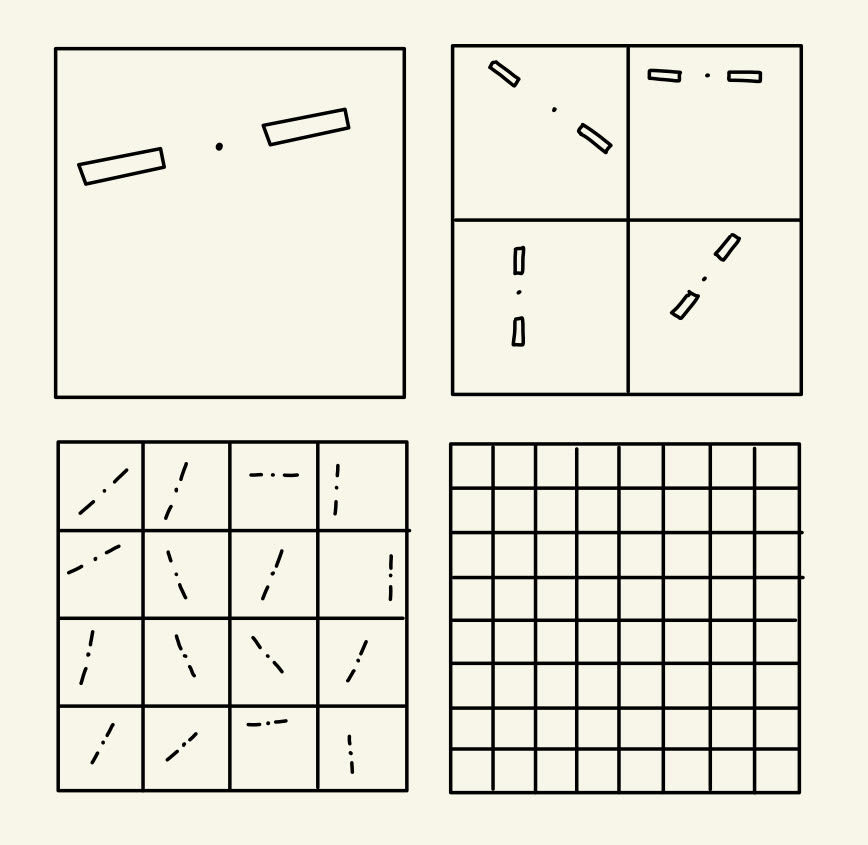}
		\caption{Definition of a backtrack and a backtrack for each dyadic square}
	\end{center}
\end{figure}

\section{Constructing a spiral chain}\label{Sec:spiralChainConstruction}

Fix a square $Q$ which contains no backtracks.
By translating and scaling, without loss of generality this square is $[0,1]^2$.
Draw $M$ many \textbf{radial rays} $r_1,\dots , r_M$ 
starting at the center of the square $(1/2,1/2)$ 
and having angles $\dfrac{2\pi}{M}, 2\dfrac{2\pi}{M}, \dots , 2\pi$.
The following basic observation is crucial in our construction.
(It is uniform smoothness in disguise.)

\textbf{Observation}: Suppose that $q$ lies on one of the rays, say $r_j$ and $L$ is the line through $q$ which is perpendicular to $r_j$.
Then the points $\{a\} = L \cap r_{j+1}$ and $\{b\} = L \cap r_{j-1}$ are further away from $(1/2,1/2)$ than $q$ by a multiplicative factor of $1+\Theta(1/M^2)$:
$$d\left(a,\left(\dfrac{1}{2},\dfrac{1}{2}\right)\right) = d\left(b,\left(\dfrac{1}{2},\dfrac{1}{2}\right)\right) = d\left(q,\left(\dfrac{1}{2},\dfrac{1}{2}\right)\right) \sec\left(\dfrac{2\pi}{M}\right)
$$
$$=d\left(q,\left(\dfrac{1}{2},\dfrac{1}{2}\right)\right) 
\left(1 + \dfrac{2\pi^2}{M^2} + o\left(\dfrac{1}{M^2}\right)\right).$$

\begin{figure}[h]
	\begin{center}
		\includegraphics[scale=0.15]{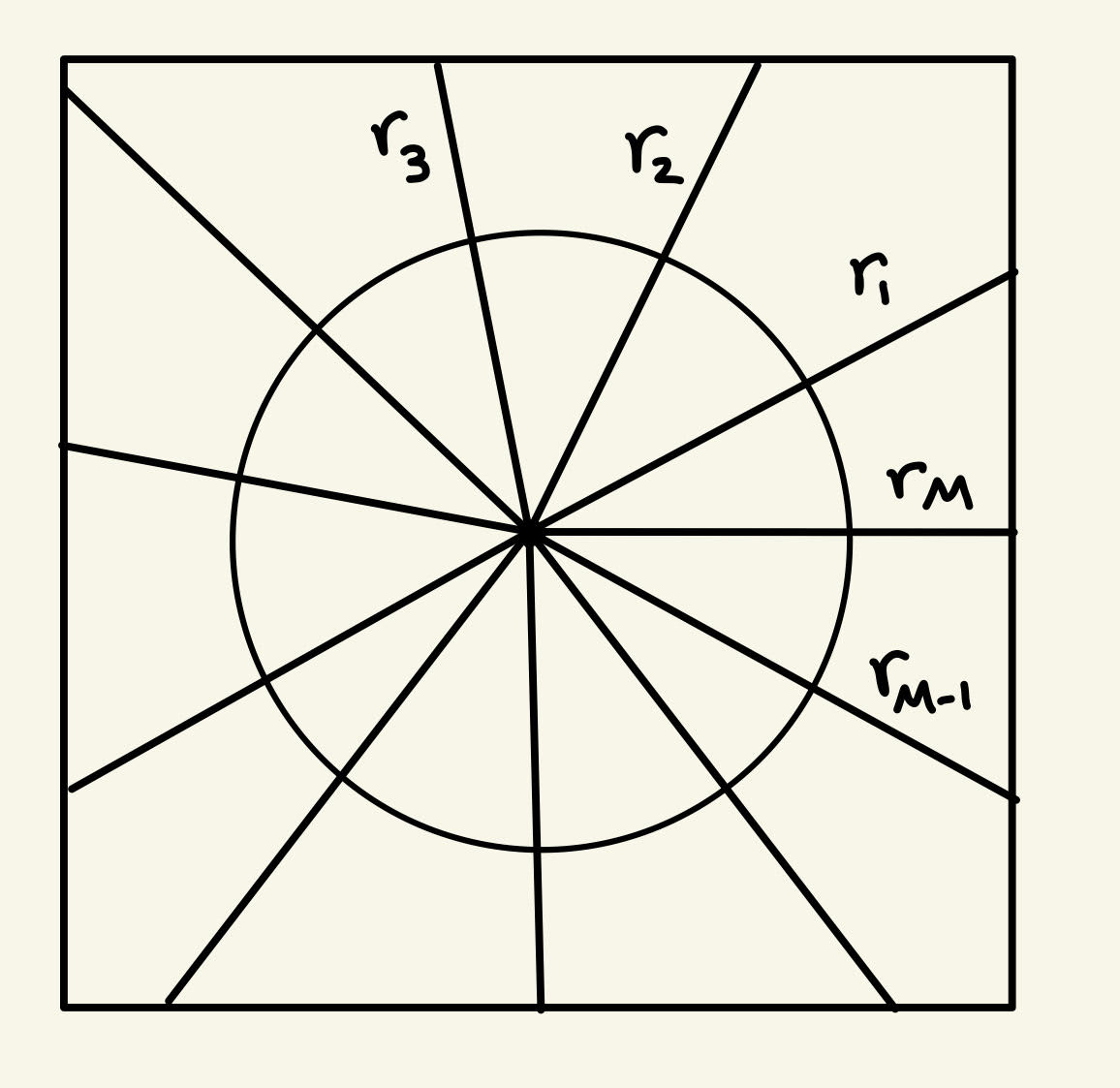}
  	\includegraphics[scale=0.15]{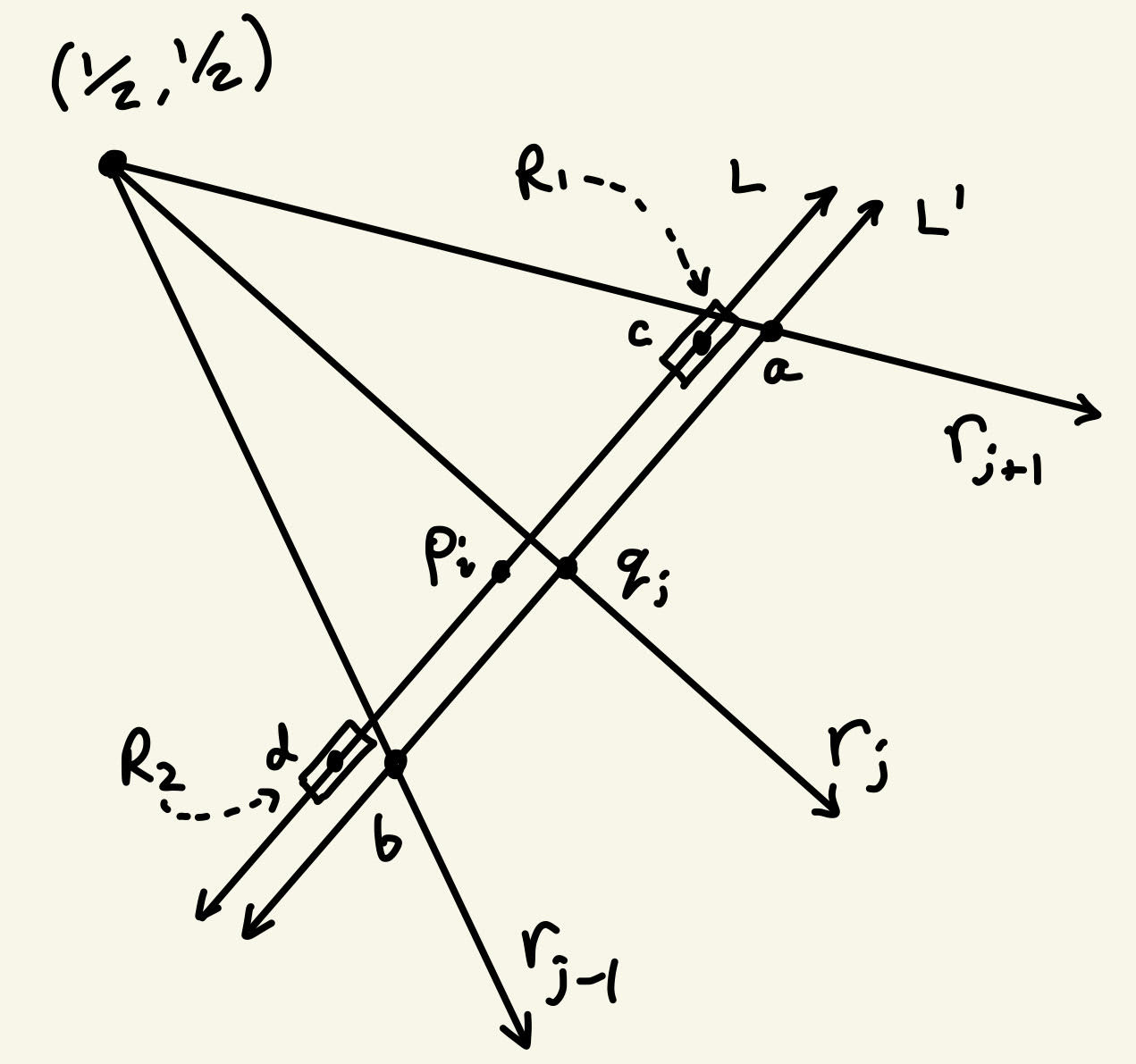}
		\caption{The radial rays and construction of the spiral chain}
	\end{center}
\end{figure}

In the next lemma, we construct a chain in the linear order $\leq$ with large jumps slowly exiting the square, which resembles a spiral.

\begin{lemma}[Spiral chain construction]\label{spiral:chain:construction}
    .\\Assume that $[0,1]^2$ has no backtrack  of length $l$ and width $w$.
    \\Suppose that $l \leq 0.01 M^{-4}$ and $M \in \n$ is divisible by $4$.
    \\\textbf{Then} there exist 
    $$p_1 > p_2 > p_3 > \dots  > p_{M^2}$$
    such that for each $i \in [M^2]$,
    $$\dfrac{1}{4}\sec\left(\dfrac{2\pi}{M}\right)^{i-1} - 2 i l \leq 
    d\left(p_i,\left(\dfrac{1}{2},\dfrac{1}{2}\right)\right)
    \leq \dfrac{1}{4}\sec\left(\dfrac{2\pi}{M}\right)^{i-1} + 2 i l,$$
    $$ \min_{j\in [M]} dist(p_i,r_j) \leq 2 i l,\;\;and$$
    $$dist(p_i,r_j) \leq 2 i l \implies \min\{dist(p_{i+1},r_{j+1}), dist(p_{i+1},r_{j-1})\} \leq 2(i+1)l.$$
\end{lemma}

\textbf{Note 1:} In words, the above 3 statements say that for each $p_i$ the distance to the center of the square is roughly $1/4(1 + 1/M^2)^i \approx 1/4 + i/4M^2$, $p_i$ is very close to one of the rays $r_1,\dots , r_M$, and the next point in the sequence $p_{i+1}$ is very close to a neighboring ray.

\textbf{Note 2:} We assume $l \leq 0.01 M^{-4}$ so that in the first inequality above, when $i = M^2$, the error is less than the increment: $2 i l < 0.02/M^2$.

\begin{proof}
Start from the point $p_1 = (3/4, 1/2)$ which lies on the radial ray $r_M$, and inductively construct the next point.
We strengthen the inductive assumption:
along with the points $p_1, p_2, p_3, \dots , p_{M^2}$ we also describe how to obtain "companion" points $q_1, q_2, q_3, \dots  , q_{M^2}$ such that for each $i$:
\\A. there is some $j \in [M]$ with $q_i \in r_j$
\\B. $d(p_i,q_i) \leq 2il$
\\C. $d\left(q_i,\left(\dfrac{1}{2},\dfrac{1}{2}\right)\right) = \dfrac{1}{4}\sec\left(\dfrac{2\pi}{M}\right)^{i-1}$
\\D. $q_i \in r_j \implies q_{i+1} \in r_{j+1} \cup r_{j-1}$.
\\(Picking the companion points $q_1, q_2, q_3, \dots  , q_{M^2}$ is not necessary, but it makes the proof less cumbersome.)
\\Start by setting $q_1 := p_1$.
Next, we describe a construction on how to get $p_{i+1}, q_{i+1}$ from $p_i, q_i$.

\textbf{The construction}:

Let $L$ be the line perpendicular to $r_{j}$ which passes through $p_i$ and $L'$ the line perpendicular to $r_{j}$ which passes through $q_i$.
Since $L \bot r_j$ and $M$ is divisible by $4$, $L$ has one of the $M$ specified angles in the definition of a backtrack.

Let $a$ and $b$ be the intersections of $L'$ with $r_{j+1}$ and $r_{j-1}$ respectively.
Let $c$ and $d$ be the two points on $L$ which have distance $d(a,q_i)$ and $d(b,q_i)$ from point $p_i$ such that the vectors $q_i a$ and $p_i c$ have the same direction, and the vectors $q_i b$ and $p_i d$ also have the same direction (see figure above).

Let $\sigma$ be the $w$-neighborhood of $L$ and denote by $R_1$ and $R_2$ the rectangles inside $\sigma$ having width $w$ and length $l$ and centers of mass $c$ and $d$ respectively.

Because we have no backtrack, there exists $p_{i+1}$ in $R_1$ or $R_2$ such that $p_i > p_{i+1}$.
If $p_{i+1} \in R_1$, put $q_{i+1}:=a$.
If $p_{i+1} \in R_2$, put $q_{i+1}:=b$.
Without loss of generality, say $p_{i+1}\in R_1$, so $q_{i+1} = a$.
This implies that $d(q_{i+1},c) = d(p_i,q_i)$.

Conditions A and D are immediate.
Condition C follows from the observation above.
For condition B, we have
$$d(p_{i+1},q_{i+1}) \leq d(p_{i+1},c) + d(c,q_{i+1}) = d(p_{i+1},c) + d(p_i,q_i) \leq l + w + 2il \leq 2(i+1)l.$$
Finally, note that the point $p_i$ stays inside the square $[0,1]^2$
so long as 
$1/4 \sec\left(2\pi/M\right)^i < 1/2-2il.$
Looking at the Taylor series of the secant function, we see that we can take $i \asymp M^2$.
\end{proof}

To each $p_1,\dots ,p_{M^2}$ we can associate $a_1,\dots ,a_{M^2} \in \z/M$ according to the index $j$ of the radial ray $r_j$ that $p_i$ is close to.
Thanks to Lemma \ref{long:walk:lemma} in the next section, we can find sets of large competitive ratio, which \textbf{zig-zag} between three different rays.

\section{A long walk on a cycle must have a zig-zag}

The following lemma is a dichotomy for any walk of length $M^2$ on the $M$-cycle.
For any $1 \leq s \leq M^{1/3}$, either there are two points of distance $s^2$ apart such that the walk zig-zags between them $M/2s$ many times, 
or else for $s^3$ consecutive steps, the walk is confined within a set of diameter $6 s^2 + 2$.

\begin{lemma}[a long walk on a cycle must have a zig zag]\label{long:walk:lemma}
    .\\Let $s, M \in \n$ with $1 \leq s \leq M^{1/3}$ and $a_1,\dots ,a_{M^2} \in \z/M$ be such that
    $a_{j+1} - a_j \in \{+1,-1\}$ for all $j = 1,\dots ,M^2$.
    \textbf{Then either} there exists $m > \dfrac{M}{s}$ and $a\in \z/M,$ and 
    $$1 \leq i_1 < j_1 < i_2 < j_2 < \dots  < i_m < j_m \leq M^2\;\;\;such\;that$$
    $$\{a_{i_1},\dots ,a_{i_m}\} = \{a\}\;\;and\;\;\{a_{j_1},\dots ,a_{j_m}\} = \{a + s^2, a - s^2\}$$
    \textbf{or else} there exists an interval $J \subset [M^2]$ of length $< s^3$, $m \geq \dfrac{s}{7}$ and $i_1, \dots . i_{m} \in J$ such that
    $$ a_{i_1} = \dots  = a_{i_m}.$$
\end{lemma}

\textbf{Illustrating examples:}
    The sequence $a_1,\dots ,a_{M^2}$ is a path on the $M$-cycle, and we wish to distinguish between the above two scenarios.
    For the "winding walk": $a_{j+1} = a_j + 1$ for every $j$, we get the first scenario.
    For the "constant walk": $a_{j+1} = a_j + 1$ for even $j$ and $a_{j+1} = a_j - 1$ for odd $j$, we get the second scenario.
    There is also a constant walk which makes one complete revolution: $a_{j+1} = a_j + 1$ for even $j$ and $a_{j+1} = a_j - 1$ for odd $j$ with the exception that for every $M$ steps we add a $1$ two times in a row before oscillating again between $+1$ and $-1$.
    This example also corresponds to the second scenario.

\textbf{A tight example:} 
    Fix $M$ and $1\leq s \leq M^{1/3}$.
    The following example shows that the parameters in the dichotomy are tight (up to multiplicative factors).
    The walk starts at $0$ and moves up to $s^2$ in $s^2$ steps.
    Then the walk moves down to $0$ in the next $s^2$ steps.
    We repeat this oscillation between $0$ and $s^2$ for $\asymp s$ many steps and then move up to $s^2$.
    We get a sub-walk which takes $\asymp s^3$ total steps, starts at $0$ and ends at $s^2$.
    Now repeat the sub-walk starting at $s^2$ and ending at $2s^2$.
    Then repeat the sub-walk from $2s^2$ to $3s^2$ and so forth for $\asymp M^2/s^3$ total iterations.
    At the end, we get a walk which winds around the circle $\asymp M/s$ times.
    The details for checking tightness are left to the reader.

\begin{proof}
    We will fix an \textbf{oscillation window} $\Delta > 0$ which we will take $\Delta = s^2$.
    For each $j_1 \in [M^2]$ we let $j_2 := \min\{j > j_1 : a_j = a_{j_1}\}$ be the next time the path visits again the point $a_{j_1}$.
    Call the time $j_1$:
    \\a \textbf{time of no oscillation} if $j_2 = \infty$, i.e. the path never visits the point $a_{j_1}$ again.
    \\a \textbf{time of small oscillation} if the walk until the next visit stays withing the oscillation window:
    $$|a_j - a_{j_1}| \leq \Delta,\;\;for\;all\;j = j_1,\dots .,j_2.$$
    a \textbf{time of large oscillation} if the walk until the next visit exits the oscillation window:
    $$|a_j - a_{j_1}| > \Delta,\;\;for\;some\;j_1 < j < j_2.$$
    Here we denote by $|\cdot|$ the distance function on the Cayley graph $\z/M$ generated by $\{\pm 1\}$. (This should not be confused with the absolute value.)
    
    Denote the times of no oscillation by $\mathcal{N}$ and the times of large oscillation by $\mathcal{L}$.
    Note that we always have $|\mathcal{N}| \leq M$.
    We split into two cases:
    \\\textbf{CASE 1:} $|\mathcal{L}| > \dfrac{M^2 s}{\Delta} = \dfrac{M^2}{s}$
    \\Then by pigeonhole there exists $a \in \z/M$ such that 
    $$m:= |\{i \in \mathcal{L}: a_i = a\}| > \dfrac{M s}{\Delta} = \dfrac{M}{s}.$$
    Denote the elements of $\{i \in \mathcal{L}: a_i = a\}$ by $i_1 < \dots  < i_m$.
    Let $j_\mu$ be the first index after $i_\mu$ that exists the $\Delta$-neighborhood of $a$ for each $\mu = 1,\dots ,m$ ($j_\mu$ exists by definition).
    This gives us the first scenario.
    \\\textbf{CASE 2:} $|\mathcal{L}| \leq \dfrac{M^2 s}{\Delta} = \dfrac{M^2}{s}$.
    \\ Since $|\mathcal{N}| \leq M \leq M^2/s$, we have $|\mathcal{L} \cup \mathcal{N}| \leq 2 M^2/s$.
    Split $[M^2]$ into the union of $\dfrac{M^2}{\Delta s}$ many intervals, each of length $\Delta s = s^3$.
    By pigeonhole, at least one of these intervals, call it $J$ has a small number of points of no or of large oscillation:
    $$|J \cap (\mathcal{L} \cup \mathcal{N})| \leq \dfrac{2M^2/s}{M^2/s^3} = 2s^2.$$
    We claim that the path during the interval $J$ has not traveled far away:
    $$\max\{|a_{j_1} - a_{j_2}| : j_1,j_2 \in J\} \leq 6 s^2 + 2$$
    To see this, call $a := a_{\min J}$ (where $\min J$ is the smallest value/time in $J$) and suppose that there is $b \in \{a_j|j\in J\}$ such that $|a-b| > 3 s^2 + 1$.
    Without loss of generality the geodesic on the $M$-cycle from $a$ to $b$ is $a, a+1, a+2, \dots , b-1, b$.
    Let $j_{\max}$ be the first time in $J$ that $b$ is visited.
    For each of the points $c = a, a+1, \dots , a+2s^2 + 1$, find the largest $j < j_{\max}$ for which $a_j = c$.
    Then by definition $j \in \mathcal{L} \cup \mathcal{N}$, which gives us $|J \cap (\mathcal{L} \cup \mathcal{N})|>2s^2$, a contradiction.
    Therefore the radius of $\{a_j: j\in J\}$ centered at $a_{\min J}$ is at most $3 s^2 + 1$, and the diameter is at most $6s^2 + 2$.
    \\By pigeonhole there exists $a \in \{a_j: j \in J\}$ such that
    $|\{j \in J: a_j = a\}| \geq \dfrac{s^3}{7 s^2} = \dfrac{s}{7}$
    which gives us the second scenario.
\end{proof}

\section{Analyzing the competitive ratio of a zig-zag}

To each $p_1,\dots ,p_{M^2}$ we associate $a_1,\dots ,a_{M^2} \in \z/M$ according to the radial ray $r_{a_i} \in \{r_1,\dots ,r_M\}$ that each point $p_i$ is closest to.
Recall the following estimates from Lemma \ref{spiral:chain:construction}:
for each $i \in [M^2]$,
$d(p_i, p_{i+1}) \gtrsim 1/M$ and $p_i$ has distance $\lesssim 1/M^2$ to one of the $M$ possible rays $\{r_1,\dots ,r_M\}$, namely to ray $r_{a_i}$.
Furthermore, for each $i_1,i_2 \in [M^2]$, we have $d(p_{i_1}, p_{i_2}) \gtrsim |a_{i_1}-a_{i_2}|/M$.

We will use the following simple upper bound on the traveling salesman cost. 
For each $S \subset [0,1]^2$ and line segments $l_1, l_2, l_3$ of lengths $|l_1|,|l_2|,|l_3|$, we have
$$tsp(S) \leq 2|l_1| + 3|l_2| + 2|l_3| + d(l_1,l_2) + d(l_2,l_3) + \sum_{p \in S} 2d(p,l_1 \cup l_2 \cup l_3).$$
(The bound follows from the following tour: the tour follows $l_1$ which has length at most $|l_1|$ and makes a detour of length $2d(p,l_1)$ for each point $p \in S$ which is closer to $l_1$ than to $l_2$ or $l_3$.
Then the tour jumps to $l_2$ by paying an cost of $|l_1| + d(l_1,l_2) + |l_2|$ and follows $l_2$ making any detours for points $p \in S$ which are the closest to $l_2$ and finally the tour does the same thing for $l_3$.)

We apply Lemma \ref{long:walk:lemma} to $a_1,\dots ,a_{M^2} \in \z/M$ and $s = \lfloor M^{1/3} \rfloor$.
We get two potential scenarios.

\begin{figure}[h]
	\begin{center}
		\includegraphics[scale=0.25]{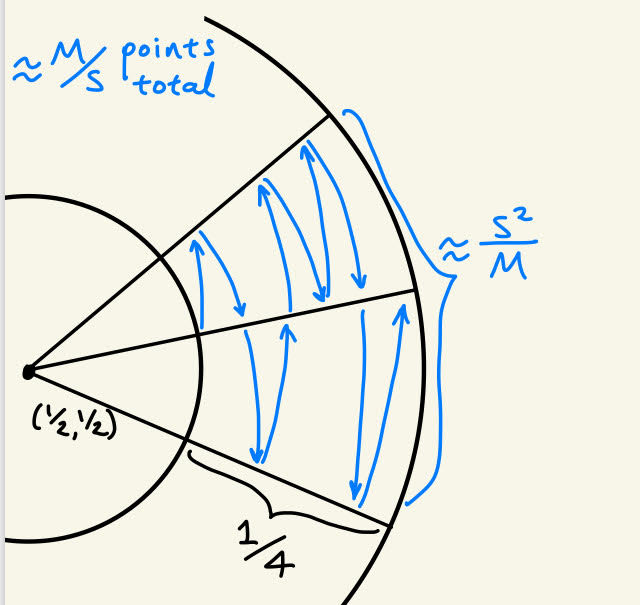}
        \includegraphics[scale=0.25]{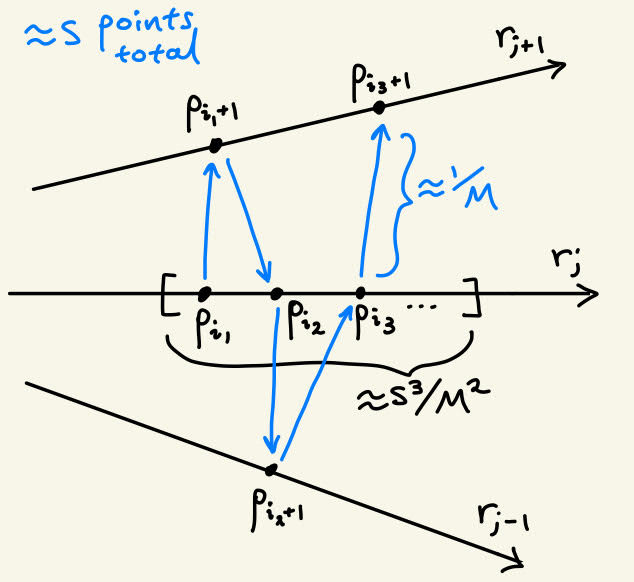}
		\caption{Applying the lemma: the first (left) and the second (right) scenarios.}
        \label{Figure-zig-zag-competitive-ratio}
	\end{center}
\end{figure}

\textbf{In the first scenario}, set $S:= \{p_{i_1},\dots ,p_{i_m},p_{j_1},\dots ,p_{j_m}\}$ and we may assume $m = \lfloor M/s\rfloor$ (else discard some of the points in $S$).
We have $d(p_{i_k},p_{j_k}) \gtrsim s^2/M$ for each $k \in [m]$ and the set $S$ is $O(1/M^2)$-close to the union three line segments of length $1/4$ (see Figure \ref{Figure-zig-zag-competitive-ratio} above). We get the bounds:
$$|S| = \dfrac{M}{s} \asymp M^{2/3},\;\;\;\;\;\;cost_{\leq}(S) \gtrsim  \dfrac{s^2}{M} \dfrac{M}{s} = s,\;\;\;\;\;\;tsp(S) \lesssim 1 + \dfrac{M}{s} \dfrac{1}{M^2} \asymp 1$$
which result in a competitive ratio
$$\dfrac{cost_{\leq}(S)}{tsp(S)} \gtrsim s \asymp \sqrt{|S|}.$$
\textbf{In the second scenario}, let $S := \{p_{i_1}, p_{i_1 + 1}, p_{i_2}, p_{i_2 +1}\dots ,p_{i_m}, p_{i_m+1}\}$ and we may assume that $m=\lfloor s/7 \rfloor$.
We have $d(p_{i_k},p_{i_k+1}) \gtrsim 1/M$ for each $k \in [m]$, and the set $S$ is $O(1/M^2)$-close to the union three line segments of length $\lesssim s^3/M^2$.
(see Figure \ref{Figure-zig-zag-competitive-ratio} above) We get the bounds:
$$|S| \asymp s,\;\;\;\;\;\;cost_{\leq}(S) \gtrsim \dfrac{1}{M} s, \;\;\;\;\;\;tsp(S) \lesssim \dfrac{s^3}{M^2} + \dfrac{1}{M^2}s \lesssim \dfrac{1}{M}$$
which result in a competitive ratio
$$\dfrac{cost_{\leq}(S)}{tsp(S)} \gtrsim \dfrac{s/M}{1/M} = s \asymp |S|.$$
We are done with \textbf{CASE B} of the proof (some dyadic square contains no backtracks).

\section{Analyzing the combination of many backtracks}

We now deal with \textbf{CASE A} of the proof.

Assume that every dyadic square has a backtrack.
Fix a backtrack for each square.
Fix a line $L$ through $[0,1]^2$ (we will choose it later at random).
For any scale $t = 0, \dots , r$, and any dyadic square $Q$ of scale $t$ with backtrack $(p,L',\sigma,R_1,R_2)$ we say that \textbf{$L$ passes through the backtrack of $Q$} if $L \cap Q$ is contained in the $\dfrac{1}{16}w 2^{-t}$ neighborhood of $L' \cap Q$ and vice versa (i.e. Hausdorff distance).

We denote by $Bad_t$ the set of all dyadic squares $Q$ of scale $t$ such that $L$ passes through the backtrack of $Q$, and call them \textbf{bad squares}.
For each scale $t$ and $Q \in Bad_{t}$, pick the point $p_Q$ corresponding to the backtrack of $Q$ (recall Definition \ref{main:definition}).

We introduce another parameter, the \textbf{scale sparsity} $c \in \n$ which we will optimize at the end to be $c \asymp \log r \asymp \log 1/w$. Define
$$ S := \bigcup_{t=0}^{\lfloor r/c \rfloor} \bigcup_{Q \in Bad_{ct}} \{p_Q\}.$$
Clearly $|S| \leq 2^{r+2}$.
We will estimate the competitive ratio of the set $S$ using the following two lemmas:

\begin{lemma}\label{tsp:bound:backtracks}
    $$tsp(S) \leq \sqrt{2} + 4 w \sum_{t=0}^{\lfloor r/c \rfloor} 
    \left[ \sum_{Q \in Bad_{ct}} \dfrac{1}{2^{ct}} \right].$$
\end{lemma}
\begin{proof}
    For any bad square $Q \in Bad_{ct}$ of scale $ct$, the point $p_Q$ of the backtrack has distance $\leq 2 w 2^{-ct}$ from the line $L$.
    Therefore, a traveling salesman tour is to follow the length of the line segment of $L \cap Q$, which is at most $\sqrt{2}$ and make a detour for each backtrack point $p$ in every bad square (moving back and forth costs at most $4 w 2^{-ct}$ for each square).
\end{proof}

\begin{figure}[h]
	\begin{center}
		\includegraphics[scale=0.26]{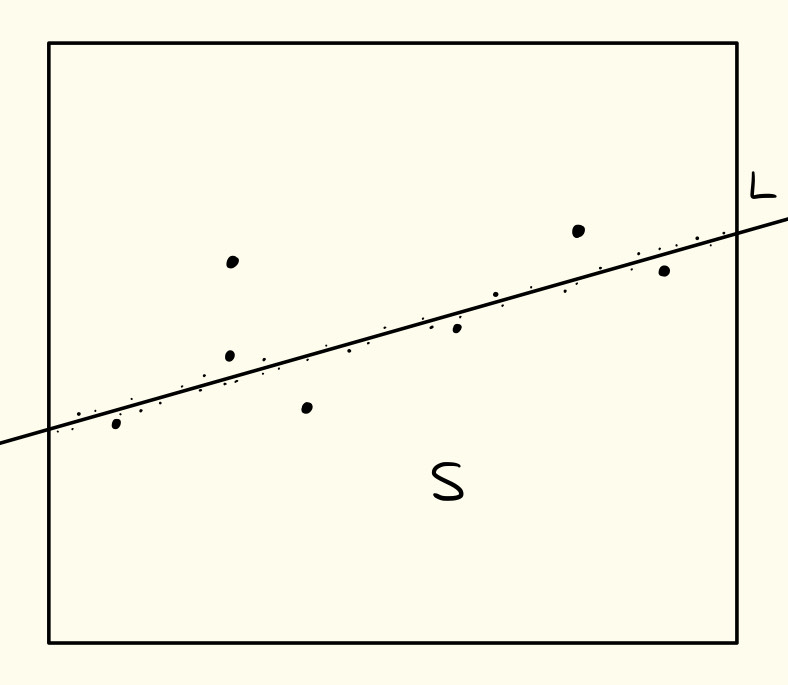}
        \includegraphics[scale=0.30]{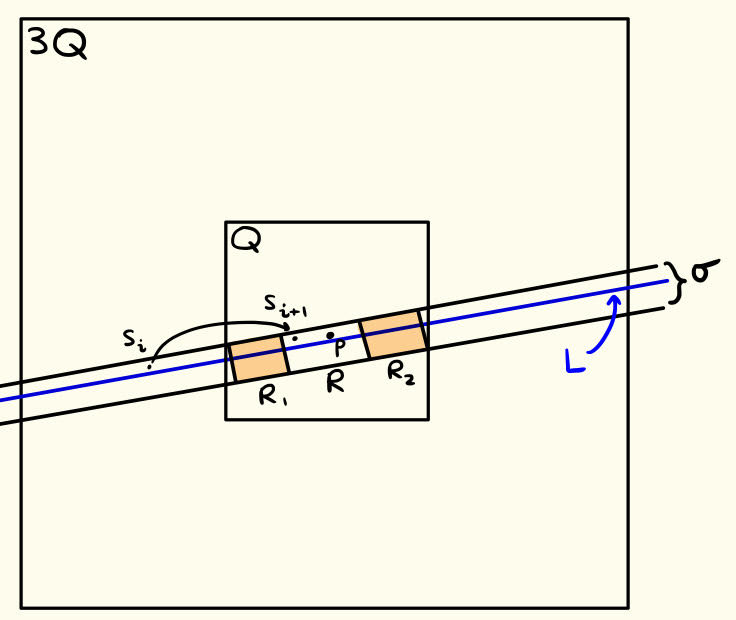}
		\caption{The backtracking set and the charging argument}
        \label{charging}
	\end{center}
\end{figure}

\begin{lemma}[Charging argument]\label{cost:bound:backtracks}
    $$2 cost_{\leq}(S) + 1 \geq l \sum_{t=0}^{\lfloor r/c \rfloor} \left( \left[\sum_{Q \in Bad_{ct}} \dfrac{1}{2^{ct}}\right] - \dfrac{18}{2^c}\right)$$
\end{lemma}
\begin{proof}
    The idea is to use a \textit{charging argument}.
    Order the points in $S$ according to $\leq$:
    $$s_1 \leq s_2 \leq \dots  \leq s_m \;\;\;\; S = \{s_1,\dots  s_m\}$$
    where $m=|S|$.
    Consider the set of all \textbf{steps} $\{s_i,s_{i+1}\}$ of the tour (where $i \in [m-1]$).
    For each $t = 0,\dots  \lfloor r/c \rfloor$ and each square $Q \in Bad_{ct}$ define its \textbf{charge} to be $l2^{-ct}$.
    For each $Q$ we will either assign the charge
    $l2^{-ct}$ to some step $\{s_i, s_{i+1}\}$ (CASE 3 below)
    or we will leave the charge of $Q$ unassigned (CASES 1 and 2 below).
    We will then argue that for each step $\{s_i, s_{i+1}\}$, the sum over all the charges assigned to it does not exceed $2 d(s_i,s_{i+1})$,
    and the sum of all the unassigned charges is small, yielding the desired lower bound.
    We have three cases:
    \\\textbf{CASE 1: If $s_1 \in Q$} 
    \\In words, the path according to $\leq$ begins from the box $Q$.
    In this case, leave the charge of $Q$ unassigned. 
    \\\textbf{CASE 2: $s_1 \notin Q$ and $3Q \cap \bigcup_{\tau=0}^{t-1} \bigcup_{Q' \in Bad_{c\tau}} \{p_{Q'}\} \neq \emptyset$.}
    \\In words, the square $3Q$ contains the point $p_{Q'}$ of a square $Q'$ larger than $Q$.
    (Here $3Q$ denotes scaling of $Q$ by $3$ which shares the same center as $Q$.)
    In this case, leave the charge of $Q$ unassigned. 
    \\\textbf{CASE 3: $s_1 \notin Q$ and $S \cap 3Q \subset \bigcup_{\tau=t}^{\lfloor r/c \rfloor} \bigcup_{Q' \in Bad_{c\tau}} \{p_Q'\}$}
    \\Let $(p,L',\sigma, R_1,R_2)$ be the backtrack of $Q$ of scale $ct$, so $\sigma$ is the strip which is the $2^{-ct} w$ neighborhood of $L'$.
    Since $3Q$ contains no point $p_{Q'}$ of scale $c\tau < ct$, all the points in $S \cap 3Q$ must be very close to the line $L$, and hence also the line $L'$. 
    The specific claim we want to make is:
    \\\textbf{Claim}: The $S \cap 3Q$ is contained in the $w2^{-ct}/2$-neighborhood of $L'$.
    \\(Let's see why this holds.
    If $L''$ is a line corresponding to the backtrack of a square $Q' \subset 3Q$ of scale $\geq ct$, then the $w2^{-ct}/16$-neighborhood of $L \cap 3Q$ contains $L''\cap Q'$ and the $w2^{-ct}/8$-neighborhood of $L' \cap 3Q$ contains $L \cap 3Q$. (We lose a factor of $2$ because we enlarged the square $Q$ to $3Q$.)
    By the triangle inequality, this means that the $w2^{-ct}/4$-neighborhood of $L'$ contains $L'' \cap Q$.
    Again by the triangle inequality, we conclude that the $w2^{-ct}/4$ neighborhood of $L'' \cap Q'$ (which contains $p_{Q'}$) is contained in the $w2^{-ct}/2$-neighborhood of $L'$, which is precisely the strip $\sigma$.)
    \\Consider the rectangular region between $R_1$ and $R_2$ and call it $R$ (which is disjoint from $R_1$ and $R_2$).
    Find the smallest $i$ such that $s_i \notin R \cup R_1 \cup R_2$ while $s_{i+1} \in R$.
    Such an $i$ exists by the definition of a backtrack; in particular $s_{i+1}$ is the smallest element in $R \cap S$ (see Figure \ref{charging}).
    
    Assign the charge of $Q$ to the step $\{s_i, s_{i+1}\}$.
    Observe that $d(s_i, s_{i+1}) \geq l 2^{-ct}$ because the point $s_i$ either belongs to the strip $\sigma$ or else it belongs outside $3Q$.
    (This is the reason we need to enlarge each square $Q$ to $3Q$.)
    
    We split the sum according to the corresponding cases
    $$\sum_{t=0}^{\lfloor r/c \rfloor} \sum_{Q \in Bad_{ct}} \dfrac{l}{2^{ct}} 
    = \sum_{t=0}^{\lfloor r/c \rfloor} \sum_{Q\in Bad_{ct}\;CASE\;1} \dfrac{l}{2^{ct}}
    + \sum_{t=0}^{\lfloor r/c \rfloor} \sum_{Q\in Bad_{ct}\;CASE\;2} \dfrac{l}{2^{ct}}$$
    $$+ \sum_{i=1}^{m-1} \sum_{t=0}^{\lfloor r/c \rfloor} \sum_{Q\in Bad_{ct}\;charges\;\{s_i, s_{i+1}\}} \dfrac{l}{2^{ct}},$$
    where by $"Q\in Bad_{ct}\;charges\;\{s_i, s_{i+1}\}"$ we mean the condition that the square $Q$ lies in CASE 3, and that the step $\{s_i, s_{i+1}\}$ is the same as the one chosen in CASE 3.
    \\\textbf{For the first sum}, for each fixed $t$, $s_1$ can only belong to at most one square of scale $t$. Thus
    $$\sum_{t=0}^{\lfloor r/c \rfloor} \sum_{Q\in Bad_{ct}\;CASE\;1} \dfrac{l}{2^{ct}} \leq \sum_{t=0}^{r-1} \dfrac{l}{2^{t+1}} \leq 2l \leq 1.$$
    \\\textbf{For the second sum}, we do not estimate it.
    For each fixed $t = 0,\dots  \lfloor r/c \rfloor$, the number of boxes in the second sum is at most 
    $$9 |\bigcup_{\tau=0}^{t-1} \bigcup_{Q \in Bad_{c\tau}} \{p_Q\}| 
    \leq 9 \sum_{\tau=0}^{t-1} 2^{c\tau} = 9\dfrac{2^{ct} - 1}{2^{c}-1} \leq 18 \dfrac{2^{ct}}{2^c}.$$
    (This holds because for each $p_Q$ of scale $<t$ we get at most $9$ boxes of scale $t$ that fall into case 2 due to $p_Q$.)
    We simply subtract the largest amount of charge we could have obtained if all these squares were in case 3.
    \\\textbf{For the third sum}, fix a step $\{s_i,s_{i+1}\}$.
    For each scale $t$, there can only be one bad square of scale $t$ which contains $s_{i+1}$ and hence there can only be one square of scale $t$ which assigns its charge to $\{s_i,s_{i+1}\}$.
    As noted above, if the charge of a square $Q$ of scale $t$ is assigned to the step $\{s_i,s_{i+1}\}$, then necessarily $d(s_i, s_{i+1}) \geq l 2^{-t}$
    (so large squares can never charge the step $\{s_i,s_{i+1}\}$).
    \\Let $t_0$ be such that $l 2^{-t_0} \leq d(s_i,s_{i+1}) < l 2^{-t_0+1}$.
    We have:
    $$\sum_{t=0}^{\lfloor r/c \rfloor} \sum_{Q\in Bad_{ct}\;charges\;\{s_i, s_{i+1}\}} \dfrac{l}{2^{ct}} \leq \sum_{t=t_0}^{r} \dfrac{l}{2^{t}} \leq \dfrac{2 l}{2^{t_0}} \leq 2 d(s_i,s_{i+1}).$$
    Summing over all steps $\{s_i, s_{i+1}\}$:
    $$ l \sum_{i=1}^{m-1} \sum_{t=0}^{r-1} \sum_{Q\in Bad_t\;charges;\{s_i, s_{i+1}\}} \dfrac{1}{2^{ct}} \leq 2 \sum_{i=1}^{m-1} d(s_i,s_{i+1})$$
    and the proof is complete.
\end{proof}

We denote the \textbf{sum of charges} by $\Sigma := \sum_{t=0}^{\lfloor r/c \rfloor} \left[\sum_{Q \in Bad_{ct}} \dfrac{1}{2^{ct}}\right]$ and combine the above Lemmas \ref{tsp:bound:backtracks} and \ref{cost:bound:backtracks} above to obtain:
$$\dfrac{cost_{\leq}(S)}{tsp(S)} \geq \dfrac{-\dfrac{1}{2} + \dfrac{l}{2} \Sigma - l \dfrac{r}{c}\dfrac{9}{2^c}}{\sqrt{2} + 4 w \Sigma}.$$
Intuitively, we want to lower bound the competitive ratio by the length-to-width ratio $l/w$.
Inspecting the above estimate, this should hold when the sum of charges $\Sigma$ is sufficiently large.
The numerics are a little tedious:
If
$$\Sigma > \dfrac{1}{w} + 90 \dfrac{r}{c 2^c}$$
then we can rewrite the inequality as
$$-\dfrac{1}{2} + \dfrac{l}{10}\Sigma \geq -\dfrac{1}{2} + \dfrac{l}{10 w} + l \dfrac{9r}{c 2^c}.$$
This inequality can be used to estimate the term:
$$-\dfrac{1}{2}  + \dfrac{l}{2}\Sigma - \dfrac{9lr}{c 2^c} \geq \dfrac{l}{10 w} + \dfrac{4 l}{10} \Sigma -\dfrac{1}{2} = \dfrac{l}{20 w}(2 + 8 w \Sigma - \dfrac{10w}{l}) > \dfrac{l}{20 w}(\sqrt{2} + 4 w \Sigma)$$
(For the last line assume that $w < 0.01 \times l$. That way $2 - 10w/l > \sqrt{2}$.)
We arrive at the competitive ratio lower bound:
$$\dfrac{cost_{\leq}(S)}{tsp(S)} \geq \dfrac{l}{20 w}.$$
provided that the sum of charges
$\Sigma = \sum_{t=0}^{\lfloor r/c \rfloor} \left[\sum_{Q \in Bad_{ct}} \dfrac{1}{2^{ct}}\right]$ 
is sufficiently large (and that $w < 0.01 \times l$).
In the final section we sample a line which satisfies this condition and optimize all the relevant constants.

\section{Taking a random line and optimizing the parameters}

We take a random line and use linearity of expectation.
Here is how we construct a random line through a $[0,1]^2$.
First, pick an angle uniformly at random from the angles $\theta \in \{\dfrac{2\pi}{M}, 2\dfrac{2\pi}{M}, \dots , 2\pi\}$.
Secondly, pick a random $y$-intercept uniformly at random among all possible $y$-intercepts which intersect $[0,1]^2$, i.e.
$$b \in \{y \sin(\theta) + x \cos(\theta) | (x,y) \in [0,1]^2\}$$
and consider the random line $L: y \sin(\theta) + x \cos(\theta) = b$.

For each dyadic square $Q$ of scale $t$ we have
$$\p(L\;passes\;through\;backtrack\;in\;Q) \geq \dfrac{w 2^{-t}}{2M}$$
(This holds because with probability $1/M$, $L$ has the same angle as the backtrack in $Q$.
Next there is probability $\geq w 2^{-t}/2$ that the random line will be $w 2^{-t}$-far away from the line of the backtrack).
\\We conclude that:
$$\e \left[\sum_{t=0}^{\lfloor r/c \rfloor} \sum_{Q \in Bad_{ct}} \dfrac{1}{2^{ct}}\right] \geq
\sum_{t=0}^{\lfloor r/c \rfloor} 2^{2ct} \dfrac{w 2^{-ct}}{2M} \dfrac{1}{2^{ct}} = \dfrac{r w}{2cM},$$
so there exists a line $L$ with $\Sigma = \sum_{t=0}^{\lfloor r/c \rfloor} \sum_{Q \in Bad_{ct}} \dfrac{1}{2^{ct}} \geq \dfrac{r w}{2cM}$.
The estimate
$$\dfrac{cost_{\leq}(S)}{tsp(S)} \geq \dfrac{l}{20 w}$$
follows as long as the parameters satisfy the two constraints:
$$\dfrac{1}{2} \dfrac{r w}{2 c M} \geq \dfrac{1}{w}\;\;\;and\;\;\; \dfrac{1}{2} \dfrac{r w}{2 c M} \geq 90 \dfrac{r}{c 2^c}.$$
The only other constraints that appear in our construction are:
$$0 < w < \dfrac{1}{100} l \leq \dfrac{M^{-4}}{10000} < 1.$$
In summary, we optimize the parameters $l>0$, $w>0$, and $c \in \n$ in the program
$$\max\{\dfrac{l}{20 w}\mid \sqrt{\dfrac{4 c M}{r}}\leq w < \dfrac{1}{100} l \leq \dfrac{1}{10000 M^4}\;\;and\;\; w \geq \dfrac{360 M}{2^c}\},$$
so we take $l = \dfrac{1}{100 M^4}$ and $w \asymp \max\{\sqrt{\dfrac{c M}{r}},\dfrac{M}{2^c}\}$.
We choose $c \in \n$ to equate (up to constant multiplicative factors) the two constraints on $w$ so:
$$\sqrt{\dfrac{c M}{r}} \asymp \dfrac{M}{2^c} \implies c \asymp \log r.$$
In particular, choosing 
$$c = \dfrac{1}{2}\log_2 r + \dfrac{1}{2} \log_2 M - \log_2 360 < \log r \;\;\; and \;\;\;w = \sqrt{\dfrac{4 M \log r}{r}},$$
both constraints of $w$ are satisfied (using the assumption on $M$ in Proposition \ref{main-prop}),
and we arrive at the promised bound:
$$\dfrac{cost_{\leq}(S)}{tsp(S)} \geq \dfrac{1}{20} \dfrac{1}{100M^4} \sqrt{\dfrac{r}{4 M \log r}}
\geq 10^{-4} \sqrt{\dfrac{r}{M^{9} \log r}}.$$
The proof of Proposition \ref{main-prop} is complete, and Theorem \ref{main:theorem} follows.

\section{Remarks}

\subsection{Exponent $1/2$ is the barrier in the proof of Theorem \ref{main:theorem}}\label{Remark:Barrier_of_1_over_2}
We give a heuristic argument that the proof and definitions for Theorem \ref{main:theorem} cannot give a lower bound on the competitive ratio which is larger than $\sqrt{\log |S|}$.
Here is the key bottleneck: the best lower bound on the sum of charges the proof could possibly give is
$$\Sigma := \sum_{t=0}^{\lfloor r/c \rfloor} \left[\sum_{Q \in Bad_{ct}} \dfrac{1}{2^{ct}}\right] \gtrsim w r.$$
This is because we are treating each scale separately and taking a random line. 
A random line cannot pass through a backtrack of scale $ct$ with probability $\gg w 2^{-ct}$.
Also, the backtracks are not correlated to each other; they are free to take any $y$-intercept.
Consider the adversarial case where all the backtracks at all scales are parallel and have a strip which is chosen independently and uniformly at random among all strips of width $w$ contained in each square.
In this case, it is impossible to pick a line $L$ (in fact, even a rectifiable subset) for which the total sum of charges is $\gg wr$.
\\(Moreover, even if we adjust the dichotomy to obtain many backtracks for each square, the argument in Lemma \ref{spiral:chain:construction} can give at most $M^{2}$ backtracks, and one can check (by the same estimate as below) that this adjustment does not help.)

Next, assuming that $\Sigma \lesssim wr$, we show that the lower bound cannot be improved.
We have 2 cases:

\textbf{CASE 1}: $w \Sigma \gtrsim 1$.
\\In this case $tsp(S) \lesssim 1 + w \Sigma \asymp w \Sigma$ and the best lower bound we get from Lemmas \ref{tsp:bound:backtracks} and \ref{cost:bound:backtracks} is at best
$$\dfrac{cost_{\leq}(S)}{tsp(S)} \gtrsim \dfrac{M^{-4}}{w\Sigma}(\Sigma- \dfrac{18 r}{c 2^c}) \gtrsim^{at\;best} \dfrac{M^{-4}}{w}.$$
On the other hand:
$$\dfrac{1}{w^2} \lesssim \dfrac{\Sigma}{w} \lesssim \dfrac{wr}{w} = r \implies \dfrac{1}{w} \lesssim \sqrt{r}.$$
So even if we set $M$ to be as small as we wish, we can never get a lower bound which is better than $\sqrt{r}$.

\textbf{CASE 2}: $w \Sigma \lesssim 1$.
\\In this case $tsp(S) \lesssim 1 + w \Sigma \asymp 1$ and the best lower bound we get from Lemmas \ref{tsp:bound:backtracks} and \ref{cost:bound:backtracks} is at best
$$\dfrac{cost_{\leq}(S)}{tsp(S)} \gtrsim M^{-4}(\Sigma- \dfrac{18 r}{c 2^c}) \gtrsim^{at\;best} M^{-4} \Sigma.$$
Now suppose that we could take $\Sigma \gg \sqrt{r}$ yielding at an improved competitive ratio.
Then $\sqrt{r} \ll \Sigma \lesssim wr$ hence $1/w \ll \sqrt{r}$.
Since $w \Sigma \lesssim 1$, we get: $\sqrt{r} \ll \Sigma \lesssim 1/w \ll \sqrt{r}$, a contradiction.


\subsection{Comparison of the definition of a backtrack with that of \cite{hajiaghayi2006improved}}\label{Remark:comparingDefinitions}

There are two differences between our definition of a backtrack and the one in \cite{hajiaghayi2006improved}.
First of all, in the definition of \cite{hajiaghayi2006improved}, there is an additional scenario where there are two points $p_1$ and $p_2$ that precede an entire rectangle $R$ which lies between $p_1$ and $p_2$, all near some line. 
Here we only deal with two rectangles $R_1$ and $R_2$ that succeed a single point $p$ which lies between $R_1$ and $R_2$, all near some line.
This means that the zig-zag argument of \cite{hajiaghayi2006improved} does not work with our new definition, and the new argument via a spiral chain is needed.
The second difference is that, because we use the spiral chain construction, \textbf{we can discretize the plane in a more flexible way}.
In the proof of \cite{hajiaghayi2006improved}, the square is split into a $k\times k$ grid and the smallest element is picked from each of the $k^2$ subsquares. 
This discretization helps with their zig-zag argument, but gives bad estimates in the backtracking set analysis.
Here the discretization is different: the lines must have one of $M$ different angles, and $M$ grows much slower than the length-to-width ratio $l/w$ of the backtrack.
There is also a scale sparsity parameter $c$, which ensures the backtracking set analysis works, and costs us the doubly logarithmic factor. (\cite{hajiaghayi2006improved} do not need to choose a sparsity parameter because they already work with a sufficiently high scale base $k$.)
So we have split the role of the parameter $k$ in \cite{hajiaghayi2006improved} into 3 new parameters: the number of angles $M$, the length-to-width ratio $l/w$ and the scale sparsity $c$.
(At the parameter optimization in the very end, $c$ is carefully chosen and leads to the log-log factor, choosing the ratio $w/l$ leads to the root-log factor whereas $M$ is freely chosen to be as slowly growing as we wish.)

\section*{Acknowledgments.}
I thank Assaf Naor for suggesting this problem to me, for insightful discussions, guidance, and encouragement.
I thank Manor Mendel for context around this work.
Finally, I thank the anonymous referees for their detailed feedback.

\newpage
\bibliographystyle{alpha}
\bibliography{references}

\end{document}